\documentclass{amsart}
\usepackage{graphicx}
\usepackage{latexsym}
\usepackage{amsfonts}
\usepackage[all]{xy}
\usepackage[notref,notcite]{showkeys}
\usepackage{amssymb, mathrsfs, amsfonts, amsmath}
\usepackage{amsbsy}

\newtheorem{theorem}{Theorem}[section]
\newtheorem{corollary}[theorem]{Corollary}
\newtheorem{lemma}[theorem]{Lemma}

\newtheorem{remark}[theorem]{\bf Remark}
\newtheorem{definition}[theorem]{\bf Definition}

\newtheorem{observation}[theorem]{\bf Observation}
\newtheorem{question}[theorem]{\bf Question}

\newcommand{\hess}{\textrm{Hess}}
\newcommand{\grad}{\textrm{grad}}
\newcommand{\infimo}{\textrm{inf}}

\begin{document}

\title{On submanifolds with tamed second fundamental form}
\author{G. Pacelli Bessa \and M. Silvana Costa  }
\date{\today}
\maketitle

\begin{abstract}
{\noindent Based on the ideas of Bessa-Jorge-Montenegro \cite{bessa} we show that a complete submanifold $M$ with tamed second fundamental form in a complete Riemannian manifold $N$ with  sectional curvature $K_{N}\leq \kappa \leq 0$ are proper, (compact if $N$ is compact). In addition, if $N$ is Hadamard then $M$ has finite topology. We also show that the fundamental tone  is an obstruction for a Riemannian manifold to be realized as submanifold with tamed second fundamental form of a Hadamard manifold with sectional curvature bounded below. }
\end{abstract}
\section{Introduction}
\hspace{.5cm}Let $\varphi: M
\hookrightarrow N$ be an isometric immersion of a complete Riemannian
$m$-manifold $M$ into a complete Riemannian   $n$-manifold
$N$  with sectional curvature $K_{N}\leq \kappa \leq
0$. Fix a point $x_0 \in M$ and let $\rho_{M} (x) = {\rm
dist}_{M}(x_0, x)$ be the  distance function on $M$ to $x_0$. Let
$\{C_{i}\}_{i=1}^{\infty}$ be an exhaustion sequence of $M$ by
compacts sets with $x_0 \in C_0$. Let $\{a_{i}\}\subset [0,\infty]$ be a non-increasing
sequence of possibly extend numbers defined  by
\[
\begin{array}{ccl}
a_i = \sup \left \{ \displaystyle
\frac{S_{\kappa}}{ C_{\kappa}}(\rho_{M} (x))\cdot
\Vert \alpha (x)\Vert, \, x \in M \backslash C_i \right \}.
\end{array}
\]
Where
\begin{equation}\label{eqSk}
 S_{\kappa} (t)=\left \{
\begin{array}{ccl}
\displaystyle \frac {1}{\sqrt{-\kappa}}\sinh(\sqrt{ -\kappa}\,t),&if& \kappa<0 \\
t, &if& \kappa=0
\end{array} \right.
\end{equation}
$ C_{\kappa}(t)= S_{\kappa}'(t)$ and $\Vert \alpha(x)\Vert $ is the  norm of the second
fundamental form at $\varphi(x)$. The number  $a(M)=\displaystyle\lim_{i \to
\infty}a_i$ is independent of the exhaustion sequence $\{C_{i}\}$ nor on the base point $x_{0}$.

\begin{definition}
An immersion $\varphi: M \hookrightarrow N$ of a complete Riemannian
$m$-manifold $M$ into an  $n$-manifold $N$ with sectional
curvature $K_{N} \leq \kappa \leq 0$ has tamed second fundamental
form if $a(M)< 1$.
\end{definition}
In \cite{bessa}, Bessa, Jorge and Montenegro showed
that a complete submanifold  $\varphi:M\hookrightarrow\mathbb{R}^{n}$ with tamed second fundamental
form is  proper and has finite topology, where finite topology means that $M$ is $C^{\infty}$-
diffeomorphic to a  compact smooth manifold $\overline{M}$ with boundary.  In this paper we show that Bessa-Jorge-Montenegro ideas can be adapted to show that a complete submanifold $ M\hookrightarrow N$ with tamed second fundamental form is proper. In addition if $N$ is a Hadamard manifold then $M$ has finite topology.
We prove the following theorem.

\begin{theorem}\label{logan}
Let $\varphi: M \hookrightarrow N$ be an isometric immersion of a complete
$m$-manifold $M$ into complete Riemannian $n$-manifold $N$  with sectional
curvature $K_{N}\leq \kappa \leq 0$. Suppose that $M$ has tamed
second fundamental form. Then

\begin{itemize}
\item[a.] If $N$ is compact then $M$ is compact.
\item[b.] If $N$ is noncompact then $\varphi$ is proper.
\item[c.] If $N$ is a Hadamard manifold then $M$ has finite topology.

\end{itemize}
\end{theorem}

Our second result shows  that the fundamental tone  $\lambda^{\ast}(M)$ can be an  obstruction for a
Riemannian manifold $M$ to be realized as a submanifold with tamed second
fundamental form in a Hadamard manifold with bounded sectional curvature. The fundamental tone of  a Riemannian
manifold $M$ is given by
\begin{equation}
\lambda^*(M) = \infimo \left \{\frac{\smallint_{M}\vert \grad
f\vert^2}{\smallint_{M}f^2},f \in H_{0}^1(M) \backslash \{ 0\} \right \},
\end{equation}where $H_{0}^1(M)$ is the completion of
$C_0^{\infty}(M)$ with respect to the norm $
\vert f\vert^2 = \int_{M}f^2 + \int_{M}\vert \grad f\vert^2.
$
We prove the following theorem.

\begin{theorem}\label{farias brito}
Let $\varphi: M \hookrightarrow N$ be an isometric immersion  of a complete
$m$-manifold $M$ with $a(M)<1$ into a Hadamard $n$-manifold $N$   with sectional
curvature  $\mu\leq K_{N} \leq 0$. Given $c$,  $a(M)<c<1$, there exists $l=l(m, c)\in \mathbb{Z}_+$ and a
positive constant $C=C(m, c,\mu)$ such that
\begin{equation}
\lambda^*(M) \leq C\cdot \lambda^*(\mathbb{N}^l(\mu))=C\cdot
(l-1)^{2}\mu^{2}/4,
\end{equation}where  $\mathbb{N}^l(\mu)$ is the $l$-dimensional simply connected space form of sectional curvature $\mu$.
\end{theorem}

\begin{remark}Jorge and Meeks  in 
\cite{jorge-meeks} showed that complete $m$-dimensional submanifolds $M$ of
$ \mathbb{R}^{n}$ homeomorphic to a compact Riemannian manifold
$\overline{M}$ punctured at finite number of points
$\{p_{1},\ldots,p_{r}\}$ and having a well defined normal vector at
infinity have $a(M)=0$. This class of submanifold includes the complete minimal surfaces $M^{2}\hookrightarrow \mathbb{R}^{n}$ with finite total curvature $\int_{M}\vert K\vert <\infty$ studied by Chern-Osserman \cite{chern-osserman}, \cite{osserman}, the
 complete surfaces $M^{2}\hookrightarrow \mathbb{R}^{n}$ with finite
total scalar curvature  $\smallint_{M}\vert \alpha
\vert^{2} dV<\infty $ and nonpositive curvature with respect to
every normal direction studied by White \cite{white}  and the $m$-dimensional
minimal submanifolds
 $M^{m} \hookrightarrow \mathbb{R}^{n}$ with finite
total scalar curvature  $\smallint_{M}\vert \alpha
\vert^{m} dV<\infty $ studied by Anderson \cite{anderson}. As corollary of Theorem (\ref{farias brito}) we have that  $\lambda^{\ast}(M)=0$ for any submanifold $M$ mentioned in this list above.
\end{remark}
\begin{question} It is known \cite{bessa-montenegro}, \cite{bessa-silvana} that the fundamental tones of the Nadirashvilli bounded minimal surfaces \cite{Nadirashvili} and the Martin-Morales cylindrically bounded minimal surfaces \cite{pacomartin} are positive. We ask if is there a complete properly immersed (minimal) submanifold of the $\mathbb{R}^{n}$ with positive fundamental tone $\lambda^{\ast}>0$.
\end{question}

\section{Preliminaries}

Let $\varphi:M \hookrightarrow N$ be an isometric immersion, where $M$ e $N$
are complete  Riemannian manifolds. Consider a smooth function
$g:N \rightarrow \mathbb{R}$ and the composition $f
=g\circ\varphi:M\rightarrow {\mathbb{R}}$. Identifying $X $ with $d
\varphi(X)$ we have at  $q \in M $ and for every $X \in T_qM$ that
\[
\langle {\rm\grad} f, X \rangle=df(X) = dg(X)=\langle {\rm\grad} g ,
X \rangle.
\]
Hence we write
\[
{\rm\grad} g = \grad f + ({\rm\grad} g)^{\bot},
\]
where $({\rm\grad} g)^{\bot}$ is perpendicular to $T_{q}M$. Let
$\nabla$ and $\bar \nabla$ be the Riemannian connections on $M$ e
$N$ respectively, $\alpha(x)(X,Y)$ and ${\rm\hess} f(x)(X,X)$ be
respectively the second fundamental form of the immersion $\varphi$
and the Hessian of $f$ at $x$ with $X, Y\in T_{x}M$. Using the Gauss
equation we have that
\begin{equation}\label{eqBF2}
{\rm\hess} f(x)(X,Y) ={\rm\hess}g(\varphi(x))(X,Y)
+\langle{\rm\grad}g,\alpha(X,Y)\rangle_{\varphi(x)}
\end{equation}
Taking the trace in (\ref{eqBF2}), with respect to an orthonormal
basis $\{e_1, ..., e_m\}$ for $T_{x}M$, we have that
\begin{equation}\label{eqBF3}
\begin{array}{ccl}
\Delta f (x) &=& \displaystyle\sum_{i=1}^{m}{\rm\hess} f(q)(e_i,e_i)\\
\\
& = & \displaystyle\sum_{i=1}^{m}{\rm\hess}g(\varphi(x))(e_i,e_i) +\langle{\rm\grad}g,\displaystyle\sum_{i=1}^{m}\alpha(e_i,e_i)\rangle.\\
\end{array}
\end{equation}
We should mention that the formulas (\ref{eqBF2}) and (\ref{eqBF3})
first appeared in \cite{jorge}. If
$g = h\circ \rho_{N}$, where $h:\mathbb{R}\rightarrow\mathbb{R}$ is
a smooth function and $\rho_{N}$ is the distance function to a fixed point in $N$,
then the equation (\ref{eqBF2}) becomes
\begin{equation}\label{eqHess}
{\rm\hess} f(x)(X,X)
=h''(\rho_{N})\langle{\rm\grad}\rho_{N},X\rangle^2 +h'(\rho_{N})[
{\rm\hess}\rho_{N}(X,X)
+\langle{\rm\grad}\rho_{N},\alpha(X,X)\rangle]
\end{equation}
Another important tool in this paper the Hessian Comparison Theorem,
see \cite{jorge} or \cite{yau}.
\begin{theorem}[Hessian Comparison Thm.]
Let $N$ be a complete Riemannian $n$-manifold and $y_0, y  \in N$.
Let $\gamma :[0,\rho_{N} (y)] \rightarrow N$ be a minimizing
geodesic joining $y_0$ and $y$,where
$\rho_N$ is the distance function to $y_0$ on $N$. Let $K_{\gamma}$
be the sectional curvatures of $N$ along $\gamma$. Denote by $\mu =
\inf K_{\gamma}$ and $\kappa =\sup K_{\gamma}$. Then for all $X \in
T_{y}N$, $X \perp\gamma'(\rho_{N}(y))$ the Hessian of $\rho_{N} $ at
$y= \gamma(\rho_{N}(y))$, satisfies
\begin{equation}\label{eqBF6}
\displaystyle\frac{C_{\mu}}{S_{\mu}}(\rho_{N}(y))\Vert
X\Vert^2 \geq {\rm\hess}\, \rho_{N}(y)
(X,X)\geq \displaystyle\frac{C_{\kappa}}{S_{\kappa}}(\rho_{N}(y))
\Vert X\Vert^2
\end{equation}
whereas ${\rm\hess}\rho_{N}(y)(\gamma',\gamma')= 0$.
\end{theorem}
\begin{observation}If $y\in {\rm cut}_{N}(y_{0})$ the inequality (\ref{eqBF6}) has to be understood in the following sense \[\displaystyle\frac{C_{\mu}}{S_{\mu}}(\rho_{N}(y))\Vert X \Vert^{2}\geq \lim_{j\to \infty} \hess \rho_{N}(y_{j})
(X_{j},X_{j})\geq \displaystyle\frac{C_{\kappa}}{S_{\kappa}}(\rho_{N}(y))
\Vert X\Vert^2 \cdot\]
For a sequence $(y_{j},X_{j})\to (y,X)\in TN,$  $y_{j}\notin {\rm cut}_{N}(y_{0})$.
\end{observation}

\section{Proof of Theorem \ref{logan}}
\subsection{Proof of items a. and b.}
Since that $a(M)< 1$, we have that for each $a(M)< c<1$, there
is $i$ such that $a_i\in (a(M),c)$. This means that there
exists a  geodesic ball $B_{M}(r_0)\subset M$, with $C_i \subset
B_{M}(r_0)$, centered at $x_0$ with radius $r_0 > 0$ such that
\begin{equation}
\displaystyle \frac {S_{\kappa}}{C_{\kappa} }(\rho_M
(x))\cdot\Vert \alpha(x)\Vert \leq c <1,\;\;\;
 for \;\; all\;\; x \in M \backslash B_{M}(r_0).
\end{equation}To fix the notation,  let $x_{0}\in M$, $y_{0}=\varphi (x_{0})$ and $\rho_{M}(x)={\rm dist}_{M}(x_{0},x)$ and $\rho_{N}(y)={\rm dist}_{N}(y_{0},y)$.
Suppose first that $\kappa=0$. Letting $ h(t)= t^2$ we have that $f(x)
= \rho_{N} (\varphi(x))^2$. By equation (\ref{eqHess}) the Hessian
of $f$ at $x\in M $ in the direction $X$ is given by
\begin{equation}
\displaystyle{\rm\hess} f(x)(X,X) = 2\,[\rho_{N} \;
{\rm\hess}\rho_{N}(X,X) +\rho_{N} \; \langle{\rm\grad}
\rho_{N},\alpha(X,X)\rangle  +\langle {\rm\grad}
\rho_{N},X\rangle^2](y),
\end{equation}where $y=\varphi (x)$.
By the Hessian Comparison theorem, we have that
\begin{equation}
{\rm\hess} \rho_{N}(y)(X,X) \geq \displaystyle \frac{1}{
\rho_{N}(y)}\Vert X^{\bot}\Vert^2,
\end{equation}
where $\langle X^{\bot},{\rm\grad} \rho_{N} \rangle =0$. Therefore
for every $x \in M \backslash B_{M}(r_0)$,
\begin{equation}\label{eqf}
\begin{array}{ccl}
{\rm\hess} f(x)(X,X) &=& 2[\rho_{N}\,{\rm\hess}\rho_{N}(X,X) +
\langle {\rm\grad}\rho_{N},X\rangle^2\\
\\
&+& \,\rho_{N}\langle {\rm\grad}\rho_{N},\alpha(X,X)\rangle] (y)\\
\\
&\geq&  2\,[ \rho_{N} \;\displaystyle\frac {1}{\rho_{N}}\Vert
X^{\bot}\Vert^2 +
 \Vert X^{\top} \Vert^2 + \rho_{N} \langle {\rm\grad}\rho_{N},\alpha (X,X)\rangle] (y)\\
\\
&\geq& 2\,[\Vert X^{\top}\Vert^2 + \Vert X^{\bot}\Vert^2 - \rho_{M}\,\Vert \alpha\Vert\cdot\Vert X\Vert^2] \\
\\
&\geq& 2( 1 - c)\Vert X\Vert^2
\end{array}
\end{equation}In the third to the fourth line of (\ref{eqf}) we used that $\rho_{N}(\varphi (x))\leq \rho_{M}(x)$.
If $\kappa < 0$,
 we let $h(t) = \cosh(\sqrt{ -\kappa}\,t)$ then  $f(x)= \cosh (\sqrt{ -\kappa}\,
\rho_{N})(\varphi(x))$. By equation $(\ref{eqHess})$ the Hessian of
$f$ is given by
\begin{eqnarray}\hess  f(x)(X , X ) &=& \left[-\kappa \cosh (\sqrt{ -\kappa}\,\rho_{N})\langle \grad \rho_{N},X\rangle^{2} +\, \sqrt{-\kappa} \sinh (\sqrt{ -\kappa}\,\rho_{N})\hess \rho_{N}(X,X)\right.\nonumber \\
&& \nonumber \\
& +&\left.   \sqrt{-\kappa} \sinh (\sqrt{ -\kappa}\,\rho_{N})\langle \grad \rho_{N}, \alpha (X,X)\rangle\right] (\varphi (x)).\label{eqcoth}
\end{eqnarray}
By  Hessian Comparison theorem we have that
\begin{equation}
{\rm\hess}\rho_{N}(y)(X,X) \geq \sqrt{-\kappa}\,\frac{\cosh (\sqrt{-\kappa}\rho_{N})}{\sinh (\sqrt{-\kappa}
\rho_{N})} \Vert X^{\bot} \Vert^2.
\end{equation}
Since $a(M) < 1$ we have then
\begin{equation}
\begin{array}{ccl}
\Vert \alpha (x)\Vert \leq c\,\displaystyle \sqrt{-\kappa}
\frac{\cosh(\sqrt{-\kappa}\rho_{M})}{\sinh(\sqrt{-\kappa}\rho_{M})}(x) \leq c \,\sqrt{-\kappa}\frac{\cosh(\sqrt{-\kappa}
\rho_{N})}{\sinh(\sqrt{-\kappa}\rho_{N}}(\varphi(x))
\end{array}
\end{equation}
for every $x \in M \backslash B_{M}( r_0)$ and some $c\in (0, 1)$.
The last inequality follows from the fact that $\rho_{N}(\varphi(x))
\leq \rho_{M}(x)$  and that the function $\sqrt{-\kappa}\coth(\sqrt{-\kappa}\,t)$ non-increasing.
Substituting in the equation (\ref{eqcoth}), we obtain
\begin{equation}
\begin{array}{ccl}
{\rm\hess} f(x)(X, X) &\geq & -\kappa \cosh(\sqrt{-\kappa}\rho_{N}) \Vert
X^{\bot}\Vert^2
-\kappa \cosh(\sqrt{-\kappa}\rho_{N})\Vert X^{\top} \Vert^2\\
& &\\
&+& \kappa\cdot c\cdot\cosh(\sqrt{-\kappa}\rho_{N})\Vert X \Vert^2
\\
&& \\
&\geq &-\kappa \cdot\cosh( \rho_{N})(1 - c)\Vert X \Vert^2\\
\\
&\geq&-\kappa \cdot (1 - c)\cdot \Vert X \Vert^2.
\end{array}\label{eqf2}
\end{equation}
Let $\sigma : [0,\rho_{M}(x)] \rightarrow M$ be a minimal geodesic
joining $x_0$ to $x$. For all $t > r_0 $  we have  that $(f\circ
\sigma)''(t) = \hess f(\sigma(t))(\sigma' , \sigma')\geq 2(1-c)
$ if $\kappa=0$ and  $(f\circ
\sigma)''(t)\geq -\kappa (1-c)$ if $\kappa<0$.

\noindent For $t \leq r_0$ we have that $(f\circ\sigma)''(t) \geq  b
= \infimo \left\{{\rm\hess} f(x)(\nu,\nu),x \in B_{M}(r_0), \vert
\nu \vert = 1\right\}$. Hence ($\kappa=0$)
\begin{equation}
\begin{array}{ccl}
(f\circ \sigma)'(s) &=& (f\circ \sigma)'(0) +  \int_0^s (f \circ \sigma)''(\tau)d\tau\\
\\
&\geq& (f\circ \sigma)'(0)+  \int_0^{r_0} b \,d\tau + \int_{r_0}^s 2(1 - c) d\tau\\
\\
&\geq& (f\circ \sigma)'(0)+ b\,r_0 + 2(1 - c)(s - r_0).\\
\end{array}
\end{equation}
Now, $\rho_{N}(\varphi(x_0))= {\rm dist}_{N}(y_0 , y_0)=0$ then
$(f\circ \sigma)'(0) = 0$,
 and $f(x_0)=0$, therefore
\begin{equation}
\begin{array}{ccl}
f(x)&=&  \int_0^{\rho_{M}(x)} (f\circ \sigma)'(s)ds\\
\\
&\geq&  \int_0^{\rho_{M}(x)}\left \{  b\,r_0 + 2(1 -c)(s - r_0)\right\}ds\\
\\
&\geq&  b\,r_0 \,\rho_{M} (x) + 2(1-c)(\displaystyle\frac{\rho_{M}^{2}(x)}{2} - r_0\,\rho_{M}(x))\\
\\
&\geq& (1 - c)\, \rho_{M}^{2}(x) + ( b - 2(1-c))\,r_0\,\rho_{M} (x) \\
\end{array}
\end{equation}
Thus \begin{equation} \rho_{N}^2( \varphi(x)) \geq  (1-c)\,\rho_{M}^2 (x)  +
(b - 2(1 -c))r_0\,\rho_{M} (x) \label{eqP1}\end{equation}  for all $x \in M$.
Similarly, for $\kappa <0$ we obtain that
\begin{equation}
\cosh ( \sqrt{-\kappa}\,\rho_{N})( \varphi(x)) \geq  \sqrt{-\kappa}(1 - c) \rho_{M}^2 (x)  +
(b/\sqrt{-\kappa} -\sqrt{-\kappa}(1-c))r_0\rho_{M} (x) + 1.\label{eqP2}
\end{equation} If $N$ is compact (bounded) the  righthand side of(\ref{eqP1}) and (\ref{eqP2}) inequalities are bounded above. That implies that $M$ must be compact. In fact, we can find $\mu=\mu ({\rm diam}(N),c, \kappa)$ so that ${\rm diam}(M)\leq \mu $.   Otherwise (if $N$ is complete noncompact) then if $\rho_{M}(x)\to \infty$,  then $\rho_{N}(\varphi)\to \infty $ and $\varphi $ is proper.

\subsection{Proof of item c.}
Recall that we have by hypothesis that $\varphi : M \hookrightarrow N $ is  a  complete $m$-dimensional
submanifold  with tamed second fundamental form immersed in complete $n$-dimensional Hadamard manifold $N$ with $K_{N}\leq \kappa\leq 0$.
 We can assume that $M$ is noncompact. Moreover, by the item a., proved in the last subsection, $\varphi$ is a proper immersion.
We can suppose that the extrinsic distance function of $M$ defined
by $R(x) = \rho_{N} (\varphi(x))$ is a Morse function on $M$. Let $ B_{N}(
r_0)$ the geodesic ball of $ N$ centered at $y_{0}$ with radius
$r_{0}$ and $S_{r_0} =
\partial B_{N}(r_0)$. Since  $\varphi$ is proper and $a(M)<1$ we
 can take $r_0$ so that
\begin{equation}
\displaystyle \frac {S_{\kappa}}{C_{\kappa} }(\rho_{M} (x))\Vert \alpha(x)\Vert \leq c <1,\;\;\;
 for \;\; all\;\; x \in M \backslash \varphi^{-1}(B_{N}(r_0))
\end{equation}
and $\Gamma_{r_0} = {\varphi(M)} { \displaystyle \cap} {S_{r_0}}
\not = \emptyset $ is  a submanifold of ${\rm dim} \,\Gamma_{r_0}=
m-1$. For each $y\in \Gamma_{r_0}$, let us denote by
$T_y\Gamma_{r_0} \subset T_y\varphi(M)$ the tangent spaces of
$\Gamma_{r_0}$ and $\varphi(M)$ at $y$, respectively. Since the
dimension   ${\rm dim}\, T_y\Gamma_{r_0} = m - 1$ and  ${\rm
dim}\,T_y\varphi(M)=m$, there exist only one unit vector $\nu(y) \in
T_y\varphi(M)$  such that $ T_y\varphi(M) = T_y\Gamma_{r_0} \oplus
[[\nu(y)]], $ with $\langle \nu(y),\grad \rho_{N}(y)\rangle > 0$.
This defines a smooth vector field $\nu$ on a neighborhood $V$
of $\varphi^{-1}(\Gamma_{r_0})$. Here $[[\nu(y)]]$ is the
vector space generated by $\nu(y)$. Consider the function on
$\varphi(V)$ defined by
\begin{equation}\label{psi}
\begin{array}{ccl}
\psi(y)=\langle \nu,\grad \rho_{N}\rangle(y) = \langle
\nu,{\rm\grad}R \rangle (y)= \nu(y)(R),\,y=\varphi (x).
\end{array}
\end{equation}
Then $\psi(y) = 0$ if and only if every $x = \varphi^{-1}(y) \in V $
is a critical  point of the extrinsic distance function $R$. Now for
each $y \in \Gamma_{r_0}$ fixed, let us consider the solution
$\xi(t,y) $  of the following Cauchy problem on $\varphi(M)$:
\begin{equation}\label{cauchy}
\left\{
\begin{array}{ccl}
\xi_t(t,y)&=&\displaystyle\frac{1}{\psi}\,\nu(\xi(t,y))\\
\\
\xi(0,y) &=& y\\
\end{array}
\right.
\end{equation}

We will prove that along of the integral curve $t \mapsto\xi(t,y)$
there are no critical points for $R= \rho_{N} \circ \varphi$. For
this, consider the function $(\psi \circ \xi)(t,y)$ and observe that
\begin{equation}
\begin{array}{ccl}
\psi_t &=& \xi_t \langle {\rm\grad} \rho_{N} , \nu \rangle \\
\\
&=& \langle \bar \nabla_{\xi_t} {\rm\grad} \rho_{N},\nu \rangle + \langle  {\rm\grad} \rho_{N}, \bar \nabla_{\xi_t}\nu\rangle\\
\\
&=&\displaystyle\frac{1}{\psi} \langle \bar \nabla_{\nu} {\rm\grad}
\rho_{N}, \nu  \rangle +
 \displaystyle \frac{1}{\psi}\langle  {\rm\grad} \rho_{N},\nabla_{\nu}\nu + \alpha(\nu,\nu)\rangle\\
\\
&=& \displaystyle\frac{1}{\psi} {\rm\hess} \rho_{N}(\nu,\nu) +
\frac{1}{\psi}\left[\langle {\rm\grad} \rho_{N},\nabla_{\nu}\nu \rangle + \langle {\rm\grad} \rho_{N}, \alpha(\nu,\nu)\rangle \right] \\
\\
&=& \displaystyle\frac{1}{\psi} \left[{\rm\hess} \rho_{N}(\nu,\nu) +
\langle\grad \rho_{N},\nabla_{\nu}\nu \rangle +
 \langle {\rm\grad}\rho_{N}, \alpha(\nu,\nu)\rangle \right].
\end{array}
\end{equation}
Thus
\begin{equation}\label{eqpsi}
\begin{array}{ccl}
\psi_t \psi &=& \hess \rho_{N}(\nu,\nu) + \langle\grad
\rho_{N},\nabla_{\nu}\nu \rangle + \langle\grad \rho_{N},
\alpha(\nu,\nu)\rangle
\end{array}
\end{equation}
Since  $ \langle \nu,\nu \rangle = 1$, we have at once that $\langle
\nabla_{\nu}\nu,\nu \rangle = 0$. As $\nabla_{\nu}\nu \in T_{x}M$,
we have that
$$
\langle{\rm\grad} \rho_{N},\nabla_{\nu}\nu \rangle
=\langle{\rm\grad} R,\nabla_{\nu}\nu \rangle.  $$ By equation
(\ref{psi}), we can write ${\rm \grad} R(x) =
\psi (\varphi(x))\cdot \nu(\varphi(x))$, since $\grad R(x)\perp T_{\varphi
(x)}\Gamma_{\rho_{N} (y)}$, ($\Gamma_{\rho_{N} (y)}=\varphi(M)
\displaystyle \cap \partial B_{N}(\rho_{N} (y))$). Then
$$
\langle{\rm\grad} \rho_{N},\nabla_{\nu}\nu \rangle
=\langle{\rm\grad} R,\nabla_{\nu}\nu \rangle=
 \psi \langle\nu,\nabla_{\nu}\nu \rangle =0 .
$$
Writing
\begin{equation}\label{eqnu}
\nu(y) = \cos\beta(y)\; \grad \rho_{N} + \sin\beta(y)\; \omega
\end{equation}
and
\begin{equation}\label{eqro}
\grad \rho_{N} (y) = \cos \beta\; \nu(y) + \sin \beta\; \nu^*
\end{equation}
where $\langle \omega,\grad  \rho_{N} \rangle = 0$ and $\langle \nu,
\nu^* \rangle = 0$, the equation (\ref{eqpsi}) becomes
\begin{equation}
\begin{array}{ccl}
\psi_t \psi = \sin^2\beta\; \hess \rho_{N}(\omega,\omega) + \sin\beta \; \langle \nu^*, \alpha(\nu,\nu) \rangle. \\
\end{array}
\end{equation}
From (\ref{eqnu}) we have that $\psi(y)= \cos \beta(y)$
\begin{equation}
\psi_t \psi=\sqrt{1 - \psi^2}\sqrt{1 - \psi^2}\hess\rho_{N}(\omega,\omega) + \sqrt{1 - \psi^2}\langle \nu^*,\alpha(\nu,\nu)\rangle.\\
\end{equation}
Hence
\begin{equation}
\begin{array}{ccl}
\displaystyle\frac{\psi_t \psi}{\sqrt{1 - \psi^2}} &=& \sqrt{1 - \psi^2}\hess\rho_{N}(\omega,\omega) + \langle \nu^*,\alpha(\nu,\nu)\rangle.\\
\end{array}
\end{equation}
Thus we arrive at the following differential equation
\begin{equation}\label{eqdif1}
\begin{array}{ccl}
-(\sqrt{1 - \psi^2})_t &=&\sqrt{1 - \psi^2} \; \hess\rho_{N}(\omega,\omega) + \langle \nu^*,\alpha(\nu,\nu)\rangle\\
\end{array}
\end{equation}
The Hessian Comparison Theorem implies that
\begin{equation}
\hess\rho_{N}(\omega,\omega) \geq \displaystyle \frac{C_{\kappa}}{S_{\kappa}}( \rho_{N}(\xi(t,y))).\\
\end{equation}
Substituting it in the equation $(\ref{eqdif1})$  obtain the
following inequality
\begin{equation}\label{eqdif2}
\begin{array}{ccl}
-(\sqrt{1 - \psi^2})_t &\geq& \sqrt{1 - \psi^2}\; \displaystyle \frac{C_{\kappa}}{S_{\kappa}}(\rho_{N}(\xi(t,y)))\;+ \langle
\nu^*,\alpha(\nu,\nu)\rangle.\\
\end{array}
\end{equation}
Denoting  by $R(t,y)$  the restriction of $R=\rho_{N}\circ \varphi$ to $\varphi^{-1}(\xi(t,y))$ we have \[R(t,y)= R(\varphi^{-1}(\xi(t,y))) = \rho_{N}(\xi(t,y))\] On the other hand
we have that
\begin{equation}
R_t = \langle{\rm \grad} R,\frac{1}{\psi}\nu\rangle = \langle \psi
\nu, \frac{1}{\psi}\nu \rangle = 1
\end{equation}
then
\begin{equation}
R(t,y)= t + r_0.
\end{equation}
Writing  $\displaystyle \frac{C_k}{S_k}(
\rho_{N}(\xi(t,y))) = \displaystyle \frac{C_{\kappa}}{S_{\kappa}}(t + r_0) $ in $(\ref{eqdif2})$ we have
\begin{equation}\label{eqdif3}
\begin{array}{ccl}
-(\sqrt{1 - \psi^2})_t \geq \sqrt{1 - \psi^2} \; \displaystyle\frac{C_{\kappa}}{S_{\kappa}}(t + r_0) + \langle \nu^*,\alpha(\nu,\nu)\rangle\\

\end{array}
\end{equation}
Multiplying $(\ref{eqdif3})$ by $ S_{\kappa}(t + r_0)$, obtain
\begin{eqnarray*}
\begin{array}{ccl}
-\left[S_{\kappa}(t + r_0)(\sqrt{1 - \psi^2})_t + C_{\kappa}(t + r_0)\sqrt{1 - \psi^2} \right]&\geq& S_k(t + r_0)\langle
\nu^*,\alpha(\nu,\nu)\rangle\\
\end{array}
\end{eqnarray*}
The last inequality  can be written as
\begin{equation}\label{eqdif4}
\begin{array}{ccl}
\left[S_{\kappa}(t + r_0)\sqrt{1 - \psi^2}\right]_t &\leq& - S_{\kappa}(t + r_0) \langle \nu^*,\alpha(\nu,\nu)\rangle\\
\end{array}
\end{equation}
Integrating  $(\ref{eqdif4})$ of $0$ to $t$ the resulting inequality
is the following
\[
S_{\kappa}(t + r_0)\sin\beta(\xi(t,y)) \leq S_{\kappa}(r_0)\sin\beta(y) + \int_0^t{-S_k(s + r_0)\langle\nu^*,\alpha(\nu,\nu)\rangle ds}
\]
Thus
\begin{equation}\label{eq37}
\sin\theta(\xi(t,y)) \leq\frac{S_{\kappa}(r_0)}{S_k(t + r_0)}\sin\beta(y) + \frac{1}{S_{\kappa}(t + r_0)}\int_0^t{S_{\kappa}(s +
r_0)(-\langle\nu^*,\alpha(\nu,\nu)\rangle) ds}
\end{equation}
Since  $a(M)< 1$, then
\[
- \langle\nu^*,\alpha(\nu,\nu)\rangle(\xi(s,y)) \leq \Vert\alpha(\xi(s,y))\Vert \leq c \frac{C_{\kappa}}{S_{\kappa}}(\rho_{M}(\xi(s,y))) \leq  c
\frac{C_{\kappa}}{S_{\kappa}}(\rho_{N}(\xi(s,y))) = c\frac{C_{\kappa}}{S_{\kappa}}(s + r_0)
\]
for every $s \geq 0$. Substituting in (\ref{eq37}), we have
\begin{equation}\label{eq39}
\begin{array}{ccl}
\sin\beta(\xi(t,y))&\leq&\displaystyle\frac{S_{\kappa}(r_0)}{S_{\kappa}(t + r_0)}\sin\beta(y) + \frac{c}{S_{\kappa}(t + r_0)}\int_0^t C_{\kappa}(s +
r_0)ds\\
\\
&=&\displaystyle\frac{S_{\kappa}(r_0)}{S_{\kappa}(t + r_0)}\sin\beta(y) + \displaystyle\frac {c}{S_{\kappa}(t + r_0)}(S_{\kappa}(t + r_0)-
S_{\kappa}(r_0))\\
\\
&=&\displaystyle\frac{S_{\kappa}(r_0)}{S_{\kappa}(t + r_0)}(\sin\beta(y) - c) +c  < 1\\
\end{array}
\end{equation}
for all $t \geq 0$. Therefore, along the integral curve $t
\mapsto\xi(t,y)$, there are no critical point for the function $R(x)
= \rho_{N} (\varphi(x))$ outside the geodesic ball $B_{N}(r_{0})$. Since R is a Morse function  the
critical points are isolated there are finitely many of then. In
particular, the submanifold has finitely many ends. This
concludes the proof of the Theorem (\ref{logan}).

\section{ Proof of Theorem \ref{farias brito}}
 The first ingredient for the proof of Theorem
\ref{farias brito} is the well known Barta's Theorem \cite{barta}
 stated here for the sake of completeness.

\begin{theorem}[Barta]Let $\Omega$ be a bounded
 open of a Riemannian manifold with piecewise  smooth
boundary.  Let $f \in C^2(\Omega)\cap C^0(\bar \Omega)$ with
$f|\Omega > 0$ and $f|\partial \Omega =0$. The first Dirichlet
eigenvalue $\lambda_{1}(\Omega)$ has the following bounds:
\begin{equation}
\begin{array}{ccl}
\displaystyle\sup_{\Omega} (- \frac {\Delta f}{f}) \geq \lambda_1(\Omega)
\geq  \displaystyle\inf_{\Omega}(\frac{-\Delta f}{f})\\
\end{array}
\end{equation}
With equality in $(4)$ if and only in $f$ is the first eigenfunction
of $\Omega$.\end{theorem}

 \noindent Let $\varphi: M \hookrightarrow  N$ be an isometric
immersion with tamed second fundamental form of a complete $m$-manifold $M$ into a Hadamard
$n$-manifold $N$   with sectional curvature $\mu\leq K_{N}\leq 0$.
Let $x_{0}\in M$, $y_{0}=\varphi
(x_{0})\in N$ and let $\rho_{N} (y)={\rm dist}_{N}(y_{0},y)$ be the
distance function on $N$ and $\rho_{N} \circ \varphi$
the extrinsic distance on $M$. By the proof of Theorem (\ref{logan})
there is an $r_{0}>0$ such that there is no critical points $x\in
M\setminus \varphi^{-1}(B_{N}(r_{0}))$ for $\rho_{N} \circ \varphi$, where
$B_{N}(r_{0})$ is the geodesic ball in $N$ centered at $y_{0}$ with
radius $r_{0}$. Let $R >r_{0}$ and let $\Omega \subset
\varphi^{-1}(B_{N}(R))$ be a connected component. Since $\varphi $
is proper we have that $\Omega $ is bounded with boundary $\partial
\Omega$ that we may suppose to be piecewise smooth. Let
$v:B_{\mathbb{N}^{l}(\mu)}(R)\to \mathbb{R}$ be a positive first
eigenfunction of the geodesic ball of radius $R$ in the
$l$-dimensional simply connected space form $\mathbb{N}^{l}(\mu)$ of
constant sectional curvature $\mu$, where $l$ is to be determined. The function $v$ is radial, i.e.
$v(x)=v(\vert x\vert)$, and satisfies the following differential
equation,
\begin{equation}v''(t)+(l-1)\,\frac{C_{\mu}}{S_{\mu}}(t)\,v'(t)
+\lambda_{1}(B_{\mathbb{N}^{l}(\mu)}(R))v(t)=0,
\,\, \forall\, t\in[0,R].\label{eqLambda-l}
\end{equation}With initial data $v(0)=1$, $v'(0)=0$. Moreover, $v'(t)<0$ for all $t\in (0,R]$. Where $S_{\mu}$ and $C_{\mu}$ are defined in (\ref{eqSk}) and $\lambda_{1}(B_{\mathbb{N}^{l}(\mu)}(R))$ is the first Dirichlet eigenvalue of the geodesic ball $B_{\mathbb{N}^{l}(\mu)}(R)\subset\mathbb{N}^{l}(\mu)$ with radius $R$.
Define
$\tilde{v}:B_{N}(R)\to \mathbb{R}$ by
$\tilde{v}(y)=v\circ\rho_{N}(y)$
 and $f:\Omega \to \mathbb{R}$ by $f(x)=\tilde{v}\circ \varphi(x)$. By  Barta's Theorem we have
 $\lambda_{1}(\Omega)\leq \sup_{\Omega}(-\triangle f/f)$. The Laplacian $\triangle f$ at a point $x\in M$ is given by
 \begin{eqnarray}\triangle_{M} f(x)& =&[\displaystyle\sum_{i=1}^{m}\hess\, \tilde{v}(e_i,e_i)
 + \langle \grad \tilde{v},\vec H\rangle](\varphi(x))\nonumber \\
 &=&\displaystyle\sum_{i=1}^{m}\left[{v}''(\rho_{N}
 )\langle \grad\rho_{N}, e_{i}\rangle^{2} +v'( \rho_{N}
 )\,\hess\,\rho_{N} (e_i,e_i)\right]+ v'( \rho
 )\langle \grad \rho_{N},\vec H\rangle\nonumber
 \end{eqnarray}Where $\hess\,\tilde{v}$ is the Hessian of $\tilde{v}$ in the metric
 of $N$ and $\{e_{i}\}_{i=1}^{m}$ is an orthonormal basis for $T_{x}M$ where
 we made the identification
 $\varphi_{\ast}e_{i}=e_{i}$. We are going to give an upper bound for $(-\triangle f/f)$ on $\varphi^{-1}(B_{N}(R))$.
Let $x\in \varphi^{-1}(B_{N}(R)) $ and choose an orthonormal basis $\{e_1,...,e_m\} $ for $T_{x}M$ such that $\{e_{2},\ldots, e_{m}\}$ are tangent to the distance sphere $\partial B_{N}(r(x))$ of radius $r(x)=\rho_{N} (\varphi (x))$ and $e_{1}=\langle e_{1}, \grad_{N} \bar\rho\rangle \grad_{N} \bar\rho + \langle e_{1}, \partial/\partial \theta\rangle  \partial/\partial \theta$. Where $\vert \partial/\partial \theta\vert=1$, $\partial/\partial \theta\perp \grad_{N} \bar\rho$. To simplify the notation set $t=\rho_{N}(\varphi(x))$, $\triangle_{M}=\triangle$. Then
\begin{eqnarray}\label{eq40}\triangle f(x)&=&\displaystyle\sum_{i=1}^{m}\left[{v}''(t)\langle \grad \rho_{N}, e_{i}\rangle^{2} +v'(t)\,\hess\,\rho_{N} (e_i,e_i)\right]+ v'( t)\langle \grad \rho_{N},\vec H\rangle\nonumber \\
&=&{v}''(t)\langle \grad \rho_{N}, e_{1}\rangle^{2}+v'(t)\langle e_{1}, \partial/\partial \theta\rangle^{2}\,\hess\,\rho_{N} (\partial/\partial \theta, \partial/\partial \theta) \\
&&+\displaystyle\sum_{i=2}^{m}v'(t)\,\hess\,\rho_{N} (e_i,e_i)+ v'( t)\langle \grad\rho_{N},\vec H\rangle\nonumber
\end{eqnarray}

Thus from (\ref{eq40})
\begin{eqnarray}\label{eq41}-\frac{\triangle f}{f}(x)&=&-\frac{v''}{v}(t)\langle \grad \rho_{N}, e_{1}\rangle^{2}-\frac{v'}{v}(t)\langle e_{1}, \partial/\partial \theta\rangle^{2}\,\hess\,\rho_{N} (\partial/\partial \theta, \partial/\partial \theta) \\
&&-\displaystyle\sum_{i=2}^{m}\frac{v'}{v}(t)\,\hess\,\rho_{N} (e_i,e_i)-\frac{ v}{v}'( t)\langle \grad\rho_{N},\vec H\rangle\nonumber
\end{eqnarray}
The equation (\ref{eqLambda-l}) is says that \[ -\frac{v''}{v}(t)=(l-1)\frac{C_{\mu}}{S_{\mu}}\frac{v'}{v}(t)+ \lambda_{1}(B_{\mathbb{N}^{l}(\mu)}(R))\]
By the Hessian Comparison Theorem and the fact $v'/v\leq 0$ we have  from equation (\ref{eq41}) the following inequality

\begin{eqnarray}\label{eq42}-\frac{\triangle f}{f}(x)&\leq &\lambda_{1}(B_{\mathbb{N}^{l}(\mu)}(R))][1-\langle e_{1},  \partial/\partial \theta \rangle^{2}] \\
&& -\frac{C_{\mu}}{S_{\mu}}(t)\frac{v'}{v}(t)\left[m-l+l\,\langle e_{1},  \partial/\partial \theta \rangle^{2}+\frac{S_{\mu}}{C_{\mu}}\Vert \vec H\Vert\right]\nonumber.
\end{eqnarray} On the other hand the mean curvature vector $\vec H$ at $\varphi(x)$ has norm \[ \Vert \vec H
\Vert(\varphi (x)) \leq \Vert \alpha\Vert(\varphi (x))\leq c\cdot (C_{\kappa}/S_{\kappa})(\rho_{M}(x))\leq c\cdot (C_{\kappa}/S_{\kappa})(\rho_{N}(\varphi (x))).\] We have that for any given $a(M)<c<1$ there exist $r_{0}=r_{0}(c)>0$ such that  there is no critical points $x\in M\setminus \varphi^{-1}(B_{N}(r_{0}))$ for $\rho_{N}\circ\varphi$. A critical point $x$ is such that $\langle e_{1},  \partial/\partial \theta \rangle(\varphi (x))=1 $, see equation (\ref{eqnu}), there $\langle e_{1},  \partial/\partial \theta \rangle(\varphi (x))=\sin \beta (\varphi (x))$. The inequality (\ref{eq39}) is showing that for any $x\in M\setminus \varphi^{-1}(B_{N}(r_{0}))$ we have that, ($\kappa=0$ in our case),
\begin{eqnarray}\langle e_{1},  \partial/\partial \theta \rangle(\varphi (x))&\leq &\frac{r_{0}}{\rho_{N}(\varphi (x)) + r_0}\left(\sup_{z\in\varphi^{-1}(\partial B_{N}(r_{0}))}\sin \beta(\varphi (z))) -c \right) +c \nonumber \\
&\leq & \frac{r_{0}}{r_{0} + r_0}(1 -c) +c \\
& =& \frac{1+c}{2}\nonumber
\end{eqnarray}We have then from (\ref{eq41}) the following inequality
\[ -\frac{\triangle f}{f}(x)\leq \frac{c^{2}}{4}\cdot \lambda_{1}(B_{\mathbb{N}^{l}(\mu)}(R))]
-\frac{C_{\mu}}{S_{\mu}}(t)\frac{v'}{v}(t)\left[m-l+\frac{l}{4}(1+c)^{2}+c \right]\]

\noindent
Choose the least $l\in \mathbb{Z}_{+}$ such that $ m-l+l(1+c)^{2}/4+ c\leq 0$.   With this choice of $l$ we have   for all $x\in \varphi^{-1}(B_{N}(R)\setminus B_{N}(r_{0}))$ that
\begin{equation}-\frac{\triangle f}{f}(x)\leq \frac{c^{2}}{4}\cdot\lambda_{1}(B_{\mathbb{N}^{l}(\mu)}(R)).
\end{equation}

\vspace{.2cm} \noindent Now let $x \in \varphi^{-1}(B_{N}(r_0))$. Since $1-\langle e_{1},  \partial/\partial \theta \rangle^{2}\leq 1$ and $-l+l\,\langle e_{1},  \partial/\partial \theta \rangle^{2}\leq 0$ we obtain from  (\ref{eq42}) the following inequality ($t=\rho_{N}(\varphi(x))$)
\begin{equation}\label{eq46}-\frac{\triangle f}{f}(x)\leq \lambda_{1}(B_{\mathbb{N}^{l}(\mu)}(R))] -\frac{C_{\mu}}{S_{\mu}}(t)\frac{v'}{v}(t)\left[m+\frac{S_{\mu}}{C_{\mu}}\Vert \vec H\Vert\right].
\end{equation}
 We need the following technical lemma.

 \begin{lemma}\label{lemma2}Let $v$ be the function satisfying (\ref{eqLambda-l}).   Then  $-v'(t)/t \leq \lambda_{1}(B_{\mathbb{N}^{l}(\mu)}(R))$ for all $t\in [0,R]$.
 \end{lemma}
 \noindent {\em Proof:} Consider the function $h:[0,R]\to \mathbb{R}$ given by
$h(t)=\lambda \cdot t + v'(t)$,
 $\lambda=\lambda_{1}(B_{\mathbb{N}^{l}(\mu)}(R))$. We know that $v(0)=1$, $v'(0)=0$ and $v'(t)\leq 0$ besides $v$ satisfies equation (\ref{eqLambda-l}). Observe that \[0=v''(t)+(l-1)v'+\lambda v\leq v'' +\lambda.\] Thus $v''\geq -\lambda$ and $h'(t)=\lambda +v''\geq 0$. Since $h(0)=0$ we have that $h(t)=\lambda t+v'(t)\geq 0$. This proves the lemma.

 Since that $v$ is a non-increasing positive function we have that $ v(t) \geq
v(r_0)$. Applying the Lemma (\ref{lemma2})  we obtain
\begin{eqnarray}\label{eq47}-\frac{\triangle f}{f}(x)&\leq& \lambda_{1}(B_{\mathbb{N}^{l}(\mu)}(R)) +\frac{t\cdot C_{\mu}(t)}{S_{\mu}(t)}(-\frac{v'(t)}{t})\cdot\frac{1}{v(r_{0})}\left[m+c\right]\\
&\leq & \lambda_{1}(B_{\mathbb{N}^{l}(\mu)}(R))\left[1 +r_{0}\frac{C_{\mu}}{S_{\mu}}(r_{0})\cdot\frac{1}{v(r_{0})}\left[m+
c\right]\right]
\end{eqnarray}

Thus for all $x\in \varphi^{-1}(B_{N}(R))$ we have that \begin{eqnarray}-(\triangle f/f)(x)& \leq&
 \max \left\{ \frac{c^2}{4}, \left[1 +r_{0}\frac{C_{\mu}}{S_{\mu}}(r_{0})\cdot\frac{1}{v(r_{0})}\left[m+
c\right]\right]\right\}\cdot\lambda_{1}(B_{\mathbb{N}^{l}(\mu)}(R))\nonumber \\ &=&
\left[1 +r_{0}\frac{C_{\mu}}{S_{\mu}}(r_{0})\cdot\frac{1}{v(r_{0})}\left[m+
c\right]\right]\cdot\lambda_{1}(B_{\mathbb{N}^{l}(\mu)}(R))\nonumber\end{eqnarray}
Then by Barta's Theorem \[\lambda_{1}(\Omega)\leq \left[1 +r_{0}\frac{C_{\mu}}{S_{\mu}}(r_{0})\cdot\frac{1}{v(r_{0})}\left[m+
c\right]\right]\cdot\lambda_{1}(B_{\mathbb{N}^{l}(\mu)}(R))\] Observe that $C=\left[1 +r_{0}\frac{C_{\mu}}{S_{\mu}}(r_{0})\cdot\frac{1}{v(r_{0})}\left[m+
c\right]\right]$ does not depend on $R$. So letting $R\to \infty$ we have that $\lambda^{\ast}(M)\leq C\lambda^{\ast}(\mathbb{N}^{l}(\mu))$.

\begin{corollary}[From the proof] Given $c$, $a(M)<c<1$ there exists $r_{0}=r_{0}(c)>0$,  $l=l(m,c)\in \mathbb{Z}_+$ and
 $C=C(m,\mu, c)>0$ such that  for any $R>r_{0}$ and $\Omega \subset \varphi^{-1}(B_{N}(R))$ a connected component, then
\[ \lambda^{\ast}(\Omega)\leq C\cdot\lambda_{1}(B_{\mathbb{N}^{l}(\mu)}(R)).\]
\end{corollary}

\end{document}